\documentclass[11pt]{article}
\usepackage{amsmath}
\usepackage{amsthm}
\usepackage{mathptmx}
\usepackage{amssymb}
 \usepackage{mathtools}
 \usepackage{pifont}
 \usepackage{authblk}
 \usepackage{url}
 \usepackage{multicol}
\usepackage{color}
\usepackage{multicol}
\usepackage{float}
\usepackage{tikz}
\usepackage{tikz-cd}
\usetikzlibrary{calc,fadings,decorations.pathreplacing}
\usetikzlibrary{arrows}
\usetikzlibrary{shapes.arrows}   
\usetikzlibrary{positioning}
\usepackage{caption}
\usepackage{subcaption}
\usepackage{hyperref}
\usepackage{mathrsfs}
\usepackage{wasysym}
\usepackage{mathtools}
\DeclareMathAlphabet{\mathcal}{OMS}{cmsy}{m}{n}

\definecolor{1}{rgb}{1,0.2,0.3}
\definecolor{2}{rgb}{0.1,0.3,0.5}
\definecolor{3}{rgb}{1,1,0}
\definecolor{4}{rgb}{255,255,255}
\def\ra{.7}
\newtheorem{theorem}{Theorem}[section]
\newtheorem{corollary}{Corollary}[theorem]
\newtheorem{lemma}[theorem]{Lemma}

\theoremstyle{definition}
\newtheorem{definition}{Definition}[section]
 \newdimen\R
\theoremstyle{remark}

\newcommand{\gt}{g_{\triangle}}
\newcommand{\ft}{f_{\triangle}}
\newcommand{\Tn}{\mathcal{T}_n}
\newcommand{\Ap}{{\mathcal A}}
\newcommand{\Bp}{{\mathcal B}}
\newcommand{\Cp}{{\mathcal C}}

\input{configdefns.txt}
\input{grays.txt}

\usepackage{setspace}%
\usepackage[rmargin=0.35in,lmargin=0.35in,tmargin=0.35in,bmargin=0.5in,paperwidth=5.5in,paperheight=8.5in]{geometry}
\begin{document}

\tikzset
{
  x=.23in,
  y=.23in,
}

\title{{\fontsize{14}{5} \selectfont \textbf{\MakeUppercase{Polyiamonds Attaining Extremal Topological Properties, Part II}}}}
\author{\normalsize
\textbf{Greg Malen}\\
Department of Mathematics, Duke University\\
e-mail address: gmalen@math.duke.edu
\and
\normalsize
\textbf{\'Erika Rold\'an}\\
Department of Mathematics, The Ohio State University\\
e-mail address: roldanroa.1@osu.edu
}

\date{}

\maketitle
\thispagestyle{empty}

\begin{abstract}
In Part II of this work we construct crystallized polyiamonds with $h$ holes for every $h\ge1$, that is polyiamonds which use the fewest possible tiles necessary to enclose $h$ holes. Furthermore, we prove that crystallized polyiamonds satisfy a set of structural conditions, and for every $h\ge 3$ there are multiple distinct crystallized polyiamonds with $h$ holes.
\end{abstract}
\section{Introduction}

Here and throughout, we assume that the reader is familiar with Part I \cite{malen2019polyiamonds} and only briefly recall some important statements and definitions. We denote the number of tiles and the number of holes in a polyiamond $A$ by $|A|$ and $h(A)$, respectively, and define the sequences,
\begin{equation}\label{nkhk}
h_k=3\binom{k}{2}=\frac{3}{2}k^2-\frac{3}{2}k, \qquad n_k=\frac{9}{2}k^2+\frac{3}{2}k-2,\\
\end{equation}
and
\begin{equation}
\label{def:gf}
\gt(h):= \min_{h(A)=h}|A|, \qquad \ft(n):= \min_{|A|=n}h(A).
\end{equation}
Given these sequences, our main result gives a full determination of the values of $\gt(h)$, expanding on Theorem 1.1 from Part I, which states that $\gt(h_k)=n_k$ for $k\ge 2$. As a corollary this also determines all values of $\ft(n)$.

\begin{figure}[H]
    \centering
    \begin{subfigure}[t]{.25\textwidth}
        \centering
        \begin{tikzpicture}
        \TOne[1.5]
        \end{tikzpicture}
        \subcaption{$T_1$}
    \end{subfigure}
    \begin{subfigure}[t]{.25\textwidth}
        \centering
        \begin{tikzpicture}
        \TTwo[1.5]
        \end{tikzpicture}
        \subcaption{$T_2$}
    \end{subfigure}
    \begin{subfigure}[t]{.25\textwidth}
        \centering
        \begin{tikzpicture}
        \SpirTwo[1.5]
        \end{tikzpicture}
        \subcaption{$T_3$}
    \end{subfigure}
    \caption{Polyiamonds with the minimum number of tiles for up to three holes.}
    \label{first3}
\end{figure}

\begin{theorem}\label{Ttremebundoexactly}
The first three values of $\gt$ are $\gt(1)=9$, $\gt(2)=14$, and $\gt(3)=19$. Then for $h_k \le h \le h_{k+1}$ for a fixed $k\ge2$, 
\begin{equation}
\label{eq:gfct}
\gt(h)=3h+3k+1+\left\lceil\frac{2h}{k}\right\rceil.
\end{equation}
\end{theorem}
It is straightforward to check that equation (\ref{eq:gfct}) yields that $\gt(h_k)=n_k$. And for $h\ge 3$ the values of $\gt(h)$ increase incrementally by either three or four, the latter only at $h$ such that $\lceil 2h/k\rceil$ jumps by 1. The definitions in (\ref{def:gf}) can be reformulated to state that $\ft(n)= \max\{h:\gt(h)\le n\}$, and thus Theorem \ref{Ttremebundoexactly} completely determines $\ft(n)$.

\begin{theorem}\label{ftsteps}
For all $n\ge1$, $\ft(n)=h$ where $\gt(h)\le n < \gt(h+1)$.
\end{theorem}

As a direct result of the property that $\gt(h+1)-\gt(h) \ge 1$, we see that $\ft(n+1)-\ft(n)\le 1$. While this follows incidentally here, it is not at all trivial. It may well be, \textit{a priori}, that adding a single tile to an optimal configuration allows for rearrangements which produce several new holes at once. A direct geometric proof that this cannot happen in polyiamonds can be obtained via similar methods to those used in \cite{kahle2018polyominoes} with polyominoes, in which the tiles are unit squares instead of triangles.

\begin{definition}\label{def:cris}
A polyiamond $A$ is \emph{crystallized} if $\gt(h(A))=|A|$.
\end{definition}

The spiral polyiamonds $Spir_k$ constructed in Part I \cite{malen2019polyiamonds} are crystallized, for example. A natural question to ask is whether crystallized polyiamonds are unique. In almost all cases, the answer is no. However, as was shown for $Spir_k$, their structure is rather concisely determined for all $h$.

\begin{theorem}\label{distinctcrystals}
For every positive integer $h\ge 3$ there exists at least two distinct crystallized polyiamonds with $h$ holes. 
\end{theorem}

\begin{theorem}\label{crystalstructure}
A polyiamond $A$ is crystallized if and only if the dual graph of $A$ is a tree, the outer perimeter of $A$ coincides with that of a shape of minimum perimeter for $|A|+h(A)$ many tiles, and every hole in $A$ has an area of 1.  
\end{theorem}

\newpage
\vspace*{.1cm}

\begin{figure}[H]
\centering

    \begin{subfigure}[t]{.2\linewidth}
        \centering
        \begin{tikzpicture}
            \SpirTwo[1.1]
            \FractalLayer[1.1]{3}
        \end{tikzpicture}
    \caption{$T_{9}$}
    \end{subfigure}
\hspace{.7cm}
    \begin{subfigure}[t]{.3\linewidth}
        \centering
        \begin{tikzpicture}
            \SpirTwo[1.1]
            \FractalLayers[1.1]{3}{4}
        \end{tikzpicture}
    \caption{$T_{18}$}
    \end{subfigure}
\hspace{.3cm}
    \begin{subfigure}[t]{.35\linewidth}
        \centering
        \begin{tikzpicture}
            \SpirTwo[1.1]
            \FractalLayers[1.1]{3}{5}
        \end{tikzpicture}
        \caption{$T_{30}$}
    \end{subfigure}
\\[1.5cm]
    \begin{subfigure}[t]{\linewidth}
        \centering
            \begin{tikzpicture}
                \centering
                \SpirTwo
                \FractalLayers{3}{15}
            \end{tikzpicture}
        \caption{$T_{315}$}
    \label{bigfig}
    \end{subfigure}

    \caption{Crystallized polyiamonds in the sequence $\{T_h\}_{h\ge1}$.}
\label{frac34515}
\end{figure}
\newpage
The area of a hole is defined here to be the number of tiles required to fill it. In proving these results, we iteratively construct a sequence of crystallized polyiamonds $\{T_h\}_{h\ge1}$ with $h(T_h)=h$, which satisfy the structural conditions of Theorem \ref{crystalstructure} and are distinct from $Spir_k$. The first three in this sequence are depicted in Figure 1 in Part 1 \cite{malen2019polyiamonds}, and for 9, 18, 30, and 315 holes in Figure \ref{frac34515}. 


We remind the reader of some further technical notation defined in Part I \cite{malen2019polyiamonds}. The polyiamond Hex$_k$ refers to the regular hexagon of sidelength $k$, and the $k$-th hexagonal layer is $L_k=\text{Hex}_{k}-\text{Hex}_{k-1}$. The number of interior edges of a polyiamond $A$ is $b(A)$, the number of edges bounding holes is $p_h(A)$, and the number of remaining edges on the outer perimeter is $p_{out}(A)$. The shape created by placing $n$ tiles in a hexagonal spiral gives the minimum perimeter for a polyiamond with $n$ tiles, $p_{min}(n)$.

\section{Proof of Theorem \ref{Ttremebundoexactly}}\label{FractalSequence}

\subsection{Constructing Crystallized Polyiamonds}\label{buildTh}

Starting with the central configuration $T_3$, we construct the remainder of the sequence by successively adding one of the three building blocks $\Ap$, $\Bp$, and $\Cp$ shown in Figure \ref{abcs}, with appropriate rotations, clockwise around the center.
\begin{figure}[H]
\centering
    \begin{subfigure}[t]{.25\textwidth}
        \centering
        \begin{tikzpicture}
            \coordinate (v0) at (0,0);
            \foreach \i in {1,2,3,4,5} {
            	\coordinate (v\i) at (120+60*\i:2*\ra);}
            \coordinate (v6) at ($(v3)!2!(v4)$);

            \fill[black!80] (v0) -- (v1) -- (v2) -- cycle;
            \fill[black!80] (v0) -- (v2) -- (v3) -- cycle;
            \fill[black!80] (v0) -- (v3) -- (v4) -- cycle;
            \fill[black!80] (v4) -- (v5) -- (v6) -- cycle;

            \draw[black!20, thick] (v2) -- (v5);
            \draw[black!20, thick] (v3) -- (v6);

            \draw[black!20, thick] (v1) -- (v2);
            \draw[black!20, thick] (v0) -- (v3);
            \draw[black!20, thick] (v5) -- (v4);

            \draw[black!20, thick] (v5) -- (v6);
            \draw[black!20, thick] (v1) -- (v4);
            \draw[black!20, thick] (v2) -- (v3);
        \end{tikzpicture}
    \subcaption{$\Ap$}
    \end{subfigure}
    \begin{subfigure}[t]{.25\textwidth}
        \centering
        \begin{tikzpicture}
            \foreach \i in {1,2,3,4} {
            	\coordinate (v\i) at (2*\i*\ra-6*\ra,0);}
            \coordinate (a) at (120:4*\ra);
            \coordinate (u1) at ($(v1)!1/2!(a)$);
            \coordinate (u2) at ($(a)!1/2!(v3)$);
            \coordinate (u3) at ($(u1)!2!(u2)$);

            \fill[black!80] (v1) -- (v2) -- (u1) -- cycle;
            \fill[black!80] (v2) -- (u1) -- (u2) -- cycle;
            \fill[black!80] (v3) -- (u2) -- (u3) -- cycle;
            \fill[black!80] (v3) -- (v4) -- (u3) -- cycle;

            \draw[black!20, thick] (v1) -- (u1);
            \draw[black!20, thick] (v2) -- (u2);
            \draw[black!20, thick] (v3) -- (u3);

            \draw[black!20, thick] (u1) -- (v2);
            \draw[black!20, thick] (u2) -- (v3);
            \draw[black!20, thick] (u3) -- (v4);

            \draw[black!20, thick] (v1) -- (v4);
            \draw[black!20, thick] (u1) -- (u3);
        \end{tikzpicture}
    \subcaption{$\Bp$}
    \end{subfigure}
    \begin{subfigure}[t]{.25\textwidth}
        \centering
            \begin{tikzpicture}
            \foreach \i in {1,2,3,4} {
            	\coordinate (v\i) at (2*\i*\ra-6*\ra,0);}
            \coordinate (a) at (120:4*\ra);
            \coordinate (u1) at ($(v1)!1/2!(a)$);
            \coordinate (u2) at ($(a)!1/2!(v3)$);
            \coordinate (u3) at ($(u1)!2!(u2)$);

            \fill[black!80] (v2) -- (u1) -- (u2) -- cycle;
            \fill[black!80] (v3) -- (u2) -- (u3) -- cycle;
            \fill[black!80] (v3) -- (v4) -- (u3) -- cycle;

            \draw[black!20, thick] (v2) -- (u2);
            \draw[black!20, thick] (v3) -- (u3);

            \draw[black!20, thick] (u1) -- (v2);
            \draw[black!20, thick] (u2) -- (v3);
            \draw[black!20, thick] (u3) -- (v4);

            \draw[black!20, thick] (v2) -- (v4);
            \draw[black!20, thick] (u1) -- (u3);
            \end{tikzpicture}
        \subcaption{$\Cp$}
    \end{subfigure}
    \caption{$\Ap$, $\Bp$, and $\Cp$ blocks used to build crystallized polyiamonds.}
    \label{abcs}
\end{figure}
In the $(k+1)$-th layer, either $\Ap$ or $\Bp$ is used to turn each of the 6 corners, and $\Cp$ is used to extend along the current side of $\text{Hex}_k$ being covered. The following rules indicate when to use each block, as depicted in Figure \ref{addABC}.

\begin{itemize}
    \item Use $\Ap$ when only one edge of the current side of $\text{Hex}_k$ remains uncovered. $\Ap$ is always placed with its isolated tile adjacent to the open edge of the previous block ($T_{11}$ in Figure \ref{addABC}).
    \item Use $\Bp$ when the previous block finishes by covering the last edge of a side. $\Bp$ is placed with one of the two tile edges on its long side on the open edge of the previous block, extending along the next side of $\text{Hex}_k$ ($T_{13}$ in Figure \ref{addABC}).
    \item Use $\Cp$ if at least two edges of the current side of $\text{Hex}_k$ remain uncovered. Like $\Ap$, $\Cp$ is always placed with its isolated tile adjacent to the open edge of the previous block ($T_{12}$ in Figure \ref{addABC}).
\end{itemize}

\begin{figure}[t]
    \centering
    \begin{subfigure}[t]{.2\textwidth}
        \centering
        \begin{tikzpicture}
        \SpirTwo[.95]\FractalLayer[.95]{3}\Afig[.95]{1}{4}
        \end{tikzpicture}
        \subcaption{$T_{10}$}
    \end{subfigure}
\hspace{.5cm}
    \begin{subfigure}[t]{.2\textwidth}
        \centering
        \begin{tikzpicture}
        \SpirTwo[.95]\FractalLayer[.95]{3}\Afig[.95]{1}{4}
        \GAfig[.95]{2}{4}
        \end{tikzpicture}
        \subcaption{$T_{11}$}
    \end{subfigure}
\hspace{.5cm}
    \begin{subfigure}[t]{.2\textwidth}
        \centering
        \begin{tikzpicture}
        \SpirTwo[.95]\FractalLayer[.95]{3}\Afig[.95]{1}{4}\Afig[.95]{2}{4}
        \GCfig[.95]{2}{4}{1}
        \end{tikzpicture}
        \subcaption{$T_{12}$}
    \end{subfigure}
\hspace{.5cm}
    \begin{subfigure}[t]{.2\textwidth}
        \centering
        \begin{tikzpicture}
        \SpirTwo[.95]\FractalLayer[.95]{3}\Afig[.95]{1}{4}\Afig[.95]{2}{4}\Cfig[.95]{2}{4}{1}
        \GBfig[.95]{3}{4}
        \end{tikzpicture}
        \subcaption{$T_{13}$}
    \end{subfigure}
    \caption{From $T_{10}$ to $T_{13}$, first $\Ap$ is added, then $\Cp$, and then $\Bp$.}
    \label{addABC}
\end{figure}

Using $\Cp$ adds three tiles and one hole, while $\Ap$ and $\Bp$ each add a fourth tile. In the $(k+1)$-th layer, $\Ap$ and $\Bp$ will occur at multiples of $k/2$ when $k$ is even, and if $k$ is odd they alternate at intervals of $(k-1)/2$ and $(k+1)/2$. To make this precise, enumerate the holes in the $(k+1)$-th hexagonal layer with consecutive indices starting with 1, and let $w_k(l)$ denote the index of the $l$-th $\Ap$ or $\Bp$ in this layer. 
\begin{lemma}\label{cornerindices}
For fixed $k\ge2$ and $l\in\{1,2,3,4,5,6\}$, $$w_k(l)=1+\left\lfloor (l-1)\frac{k}{2}\right\rfloor.$$
\end{lemma}
\begin{proof}
We prove this assuming that the first block in each layer is $\Ap$, which we use induction to establish. Starting with $T_3$ from Figure \ref{first3} where $k=2$, the block needed to turn the corner is $\Ap$.  Then for fixed $k\ge 3$ assume the first block in the $(k+1)$-th layer is $\Ap$, and thus $w_k(1)=1$.

This initial $\Ap$ covers the first two of the $k$ boundary edges on this side of $\text{Hex}_k$, one of which actually bounds the hole that $\Ap$ creates. Any further $\Ap$'s used in this layer, however, will only cover one boundary edge on the next side. Meanwhile, $\Bp$ and $\Cp$ always cover two boundary edges, with $\Bp$'s first edge adjacent to the previous block, and thus in the interior of $L_{k+1}$ instead of bordering $\text{Hex}_k$. These characterizations are all evident in Figure \ref{addABC}, where $T_{10}$ is also constructed from $T_9$ (Figure \ref{frac34515}) by adding the initial $\Ap$ for the fourth layer. 

Following this first $\Ap$, there are then $\lfloor(k-2)/2\rfloor$ many $\Cp$'s. If $k$ is odd there is an open edge remaining and this is followed by $\Ap$, and if $k$ is even it is followed by $\Bp$. Thus $$w_k(2)=1+\left\lfloor\frac{k-2}{2}\right\rfloor+1=1+\left\lfloor \frac{k}{2}\right\rfloor.$$

On the next side, if $k$ is odd there is only one edge covered by $\Ap$, and then $(k-1)/2$ many $\Cp$'s perfectly covering the boundary on this side, followed by $\Bp$. Then for $k$ odd, $\lfloor k/2\rfloor=(k-1)/2$, and 
\begin{align*}
    w_k(3)&=1+\left\lfloor \frac{k}{2}\right\rfloor+(k-1)/2+1;\\[5pt]
    &=1+(k-1)/2+(k-1)/2+1;\\[5pt]
    &=1+k.
\end{align*}
If $k$ is even, then $\Bp$ covers the first two edges, and $(k-2)/2$ many $\Cp$'s are placed, perfectly covering the boundary of this side, and thus followed again by $\Bp$. Then for $k$ even, $\lfloor k/2\rfloor=k/2$, and 
\begin{align*}
    w_k(3)&=1+\left\lfloor \frac{k}{2}\right\rfloor+(k-2)/2+1;\\[5pt]
    &=1+\frac{k}{2}+\frac{k}{2};\\[5pt]
    &=1+k.
\end{align*}

The next side starts with $\Bp$ in both cases, which covers two edges of the new side. This is precisely what the initial $\Ap$ on the first side did, and so the above increments repeat and prove the formula for $w_k(4)$ and $w_k(5)$. Then the pattern of the first side repeats for the fifth side, proving $w_k(6)$, and it only remains to show that the next layer must also begin with $\Ap$.

\begin{figure}[b]
    \centering
    \begin{subfigure}[t]{.4\textwidth}
        \centering
        \begin{tikzpicture}
        \GSpirTwo\GFractalLayer{3}
        \FractalLayer{4}
        \fill[black!40] (300:2*\ra) --++(0:\ra) --++(240:\ra) -- cycle;
        	\draw[black!15] (300:2*\ra) --++(0:\ra) --++(240:\ra) -- (300:2*\ra);
        \end{tikzpicture}
    \end{subfigure}
    \begin{subfigure}[t]{.4\textwidth}
        \centering
        \begin{tikzpicture}
        \SpirTwo\FractalLayer{3}
        \GAfig{1}{4}
        \Afig{2}{4}\Cfig{2}{4}{1}\Bfig{3}{4}\Afig{4}{4}\Cfig{4}{4}{1}\Bfig{5}{4}\Afig{6}{4}\Cfig{6}{4}{1}
        \end{tikzpicture}
    \end{subfigure}
    \caption{$T_{18}$ shown with its $L_4$ layer perfectly covering $\text{Hex}_3$ at left, and the uncovered edge of the initial $\Ap$ in $L_4$ on the sixth side at right.}
    \label{transitionA}
\end{figure}
Observe that the sixth side in $L_{k+1}$ follows the same pattern as the second side, and so this side of $\text{Hex}_k$ is perfectly covered by the succession of $\Cp$'s. However, as this last side transitions into $L_{k+2}$, there is also an edge from the initial $\Ap$ in $L_{k+1}$ which needs to be covered, as depicted in Figure \ref{transitionA}. Therefore the next block, which is the first of $L_{k+2}$, must be $\Ap$, and hence by induction every layer starts with $\Ap$ and the equation holds for all $k\ge 2$.
\end{proof}
\begin{lemma}\label{Th:nkhk}
$|T_{h_k}|=n_k$.
\end{lemma}
\begin{proof}
We prove this by induction. For the base case $k=2$, we have $h_2=3$ and $|T_3|=19=n_2$ (Figure \ref{first3}). For fixed $k\ge3$, assume that $|T_{h_k}|=n_k$. 

In building the $(k+1)$-th layer, we add $|L_{k+1}|=12k+6$ tiles and holes combined, starting with the last two tiles in $L_k$ and ending right before the last two in $L_{k+1}$. In the construction, we add $x+6$ holes and $3(x+6)+6$ tiles, where $x$ is the number of $\Cp$ blocks. Hence $4(x+6)+6=12k+6$, so this adds $x+6=3k$ holes and $3(3k)+6=9k+6$ tiles. It is straightforward to check that $h_k+3k=h_{k+1}$, and that $n_k+9k+6=n_{k+1}.$

Thus $|T_{h_{k+1}}|=n_{k+1}$, and by induction the statement holds for all $k\ge2$.
\end{proof}

When the $h$-th hole is in the $(k+1)$-th layer, there are $h-h_k$ blocks in this layer, and we keep track of the extra tiles added by the $\Ap$ and $\Bp$ blocks by counting $\lceil (h-h_k)/(k/2)\rceil$. Consider that
\begin{align*}
    w_k(l)-1=\left\lfloor (l-1)\frac{k}{2}\right\rfloor\le (l-1)\frac{k}{2} \ \ \ \Longrightarrow \ \ \ \frac{w_k(l)-1}{\frac{k}{2}}\le l-1.
\end{align*}
So for $h-h_k < w_k(l)$, $\lceil (h-h_k)/(k/2)\rceil\le l-1$. And
\begin{align*}
    w_k(l)=1+\left\lfloor (l-1)\frac{k}{2}\right\rfloor\ge\frac{1}{2}+(l-1)\frac{k}{2} \ \ \ \Longrightarrow \ \ \  \frac{w_k(l)}{\frac{k}{2}}\ge \frac{1}{k}+l-1.
\end{align*}
Hence for $h-h_k \ge w_k(l)$, $\lceil (h-h_k)/(k/2)\rceil\ge l$. Thus $$w_k(l)\le h-h_k < w_k(l+1) \ \ \ \Longrightarrow \ \ \ \lceil (h-h_k)/(k/2)\rceil = l.$$  For such an $h$, three tiles are added for each of the $h-h_k$ holes, plus an extra 1 for each of the $l$ corner blocks that have been placed. Then for $h \ge 3$, $|T_h|$ is given by the piecewise function defined on the intervals $3\binom{k}{2} \le h \le 3\binom{k+1}{2}$ by
\begin{equation}\label{ght:equation}
\begin{split}
|T_h|&=n_k+3\left(h-h_k\right)+\left\lceil\frac{h-h_k}{\frac{k}{2}}\right\rceil\\[8pt]
&=n_k+3\left(h-\left(\frac{3}{2}k^2-\frac{3}{2}k\right)\right)+\left\lceil\frac{2\left(h-\frac{3}{2}\left(k(k-1)\right)\right)}{k}\right\rceil\\[8pt]
&=3h+3k+1+\left\lceil\frac{2h}{k}\right\rceil.\\
\end{split}
\end{equation}

\subsection{Crystallization of $T_h$}\label{g=G}
We recall the function $M(n,h)$ from Part I, defined to be
\begin{equation}
\label{eq:Mnh}
    M(n,h)=\frac{n+2-p_{min}(n+h)}{3}.
\end{equation}
By partitioning the edges of a polyiamond $A$ into exterior edges, interior edges, and hole-bounding edges, it was shown (Lemma 2.2 in \cite{malen2019polyiamonds}) that 
\begin{equation*}
    M\left(|A|,h(A)\right)\ge h(A).
\end{equation*} 
By the monotonicty in $n$ of $M(n,h)$ (Lemma 2.3 \cite{malen2019polyiamonds}), this gives a lower bound on $\gt(h)$. We then complete the proof of Theorem \ref{Ttremebundoexactly} by proving Lemmas \ref{ght:equality} and \ref{ght:optimal}, showing that equality holds in this equation with $A=T_h$, and that $|T_h|$ is minimal.

The following are easily verified by plugging into equation (\ref{eq:Mnh}), and along with $T_1$, $T_2$, and $T_3$ in Figure \ref{first3} are a proof that $\gt(1)=9$, $\gt(2)=14$, and $\gt(3)=19$:
\begin{align*}
    &M(9,1)=1, \qquad  M(14,2)=2, \qquad M(19,3)=3,\\[5pt]
    &M(8,1)=\frac{1}{3}, \qquad  M(13,2)=\frac{4}{3}, \qquad M(18,3)=\frac{7}{3}.
\end{align*}

\begin{lemma} \label{ght:equality}
$M\left(|T_h|,h\right)=h.$
\end{lemma}
\begin{proof} 
The first three cases are shown above. We then prove this by induction, from the base case $M(19,3)=3$. Assume that $M\left(|T_{h-1}|,h-1\right)=h-1$ for some $h > 3$. Then 
\begin{equation}\label{countTh1}
\begin{split}
    3(h-1)&=|T_{h-1}|+2-p_{min}\left(|T_{h-1}|+h-1\right)\\[8pt]
    \Longrightarrow \ \ \ \ \ \  3h &= \left(|T_{h-1}|+3\right)+2-p_{min}\left(|T_{h-1}|+h-1\right).\\[8pt]
\end{split}
\end{equation}
It was noted in \cite{malen2019polyiamonds} that $p_{min}(n+1) =  p_{min}(n)-1$ if the $(n+1)$-th space in a continuous hexagonal spiral of triangles is pointing outward, and $p_{min}(n+1) =  p_{min}(n)+1$ if it is pointing inwards (see Figure 13 in \cite{yangmeyer}). Thus when adding one of the $\Ap$, $\Bp$, and $\Cp$ blocks, the increments in $p_{min}$ as the individual tiles are added one at a time are
\begin{equation}
\label{eq:abcincrements}
\begin{split}
    \Ap: +1, -1, +1, &+1, -1; \qquad \Bp: +1, +1, -1, +1, -1;\\[5pt] 
    &\Cp: +1, -1, +1, -1.
\end{split}
\end{equation}
Hence when an $\Ap$ or $\Bp$ block is used, then $|T_h|=|T_{h-1}|+4$ and \[p_{min}\left(|T_h|+h\right)=p_{min}\left(|T_{h-1}|+h-1\right)+1.\]
Thus equation (\ref{countTh1}) gives
\begin{equation}\label{AorBadd}
\begin{split}
    3h &= \left(|T_{h-1}|+3\right)+2-\left(p_{min}\left(|T_h|+h\right)-1\right)\\[8pt]
    &= \left(|T_{h-1}|+4\right)+2-p_{min}\left(|T_h|+h\right)\\[8pt]
    &=|T_h|+2-p_{min}\left(|T_h|+h\right).
\end{split}
\end{equation}
And if a $\Cp$ block is used, then $|T_h|=|T_{h-1}|+3$ and \[p_{min}\left(|T_h|+h\right)=p_{min}\left(|T_{h-1}|+h-1\right).\]
So equation (\ref{countTh1}) gives
\begin{equation}\label{Cadd}
\begin{split}
    3h &= \left(|T_{h-1}|+3\right)+2-p_{min}\left(|T_h|+h\right)\\[8pt]
    &=|T_h|+2-p_{min}\left(|T_h|+h\right).
\end{split}
\end{equation}
Equations (\ref{AorBadd}) and (\ref{Cadd}) are equivalent to $M\left(|T_h|,h\right)=h$, and hence this holds for all $h$ by induction.
\end{proof}

\begin{lemma} \label{ght:optimal}
$M\left(|T_h|-1,h\right) < h$.
\end{lemma}
\begin{proof}
Noting that each of the $\Ap$, $\Bp$, and $\Cp$ blocks end with an outward facing tile, by \ref{eq:abcincrements} we have that
\begin{align*}
    p_{min}\left(|T_h|+h\right) =  p_{min}(|T_h|-1+h)-1.
\end{align*}
Applying this equality and Lemma \ref{ght:equality}, we have that
\begin{align*}
    M\left(|T_h|-1,h\right)&=\frac{\left(|T_h|-1\right)+2-p_{min}\left(|T_h|-1+h\right)}{3}\\[6pt]
    &=\frac{|T_h|+1-p_{min}\left(|T_h|+h\right)-1}{3}\\[6pt]
    &=\frac{|T_h|-p_{min}\left(|T_h|+h\right)}{3}\\[6pt]
    &=h-\frac{2}{3}.
\end{align*}
\end{proof}
\begin{proof}[Proof of Theorem \ref{Ttremebundoexactly}]
 By Corollary 2.2.1 in Part I \cite{malen2019polyiamonds}, if there exists a polyiamond with $n$ tiles and $h$ holes where $M(n,h)=h$ and $M(n-1,h) < h$, then $\gt(h)=n$. Therefore Theorem \ref{Ttremebundoexactly} follows from Lemmas \ref{ght:equality} and \ref{ght:optimal}, and equation (\ref{ght:equation}).
\end{proof}

Thus for every $h\ge 1$, the polyiamond $T_h$ is crystallized. The properties encoded in $M(n,h)$ are then used to establish Theorem \ref{crystalstructure}.

\begin{proof}[Proof of Theorem \ref{crystalstructure}]
Let $A$ be a polyiamond with $|A|=n$ and $h(A)=h$. If these three conditions are satisfied, then $b(A)=n-1$, $p_{out}(A)=p_{min}(n+h)$, $p_h(A)=3h$, and it is straightforward that $M(n,h)=h$. Suppose instead that $n=|T_h|$ and $A$ is crystallized. Then as a corollary of Theorem \ref{Ttremebundoexactly} and Lemma \ref{ght:equality}, $M\left(|A|,h(A)\right)=h(A)$. Thus $3h = 3n-2b_{min}(n)-p_{min}\left(n+h\right).$

Each of the $\Ap$, $\Bp$, and $\Cp$ pieces ends with on outward facing tile. Thus by Lemma 2.1 in Part I \cite{malen2019polyiamonds}, the values $n+h$ for $n=|T_h|$ give thresholds for $p_{min}$ such that for any $j\ge 1$, we have $p_{min}(n+h)\le p_{min}(n+h+j).$
Now suppose $A$ has at least one hole with an area of at least 2. Then for some $j\ge0$,
\begin{equation*}
    p_{out}(A)\ge p_{min}(n+h+j)\ge p_{min}(n+h).
\end{equation*}
The extra area also increases the hole perimeter, so $p_h(A) > 3h$ and thus
\begin{equation*}
\begin{split}
    3n-2b(A)-p_{out}(A) &> 3n-2b_{min}(n)-p_{min}(n+h),\\
    \Longrightarrow \hspace{2cm} b_{min}(n) &> b(A).
\end{split}
\end{equation*}
But this is a contradiction since $|A|=n$, and thus every hole has an area of one. 

If the dual graph of $A$ is not a tree, then it contains at least $n$ edges and $b(A)\geq n > n-1 = b_{min}(n)$. If $A$ does not achieve the minimum outer perimeter, then because every hole has an area of 1, $p_{out}(A) > p_{min}(n+h)$. In either case,
\begin{equation*}
    3h = p_h(A) < 3n-2b_{min}(n)-p_{min}(n+h),
\end{equation*}
which is again a contradiction.
\end{proof}

%

Finally we address the issue of uniqueness, showing that for large enough $h$ the central configuration of $T_h$ can be swapped with some $Spir_k$.
\begin{proof}[Proof of Theorem \ref{distinctcrystals}]
For $3\le h\le 8$, examples of distinct crystallized polyiamonds are given in Figure \ref{less8}. For $h\ge 9$, observe that $T_{h_k}$ has the same outer perimeter as $Spir_k$, and its central configuration is always $Spir_2$. Hence for $h \ge h_K$, the central configuration of $T_{h}$ can be swapped out with $Spir_k$ for any $3\le k \le K$ to create distinct crystallized polyiamonds (see Figure \ref{spiralswap}).
\end{proof}
\vskip.6cm
\begin{figure}[H]
    \centering
    \begin{subfigure}[t]{.3\textwidth}
        \centering
        \begin{tikzpicture}
		\foreach \w in {3,...,6}{
		\LSide{\w}{2}
		}
		\foreach \w in {1,3,4,5}{
			\fill[black!80] (0,0) --++(60*\w:\ra) --++(120+60*\w:\ra) -- cycle;
		}
		\foreach \w in {1,2}{
			\fill[black!80] (180:\ra) --++(60+60*\w:\ra) --++(180+60*\w:\ra) -- cycle;
		}
		\fill[black!80] (240:\ra) --++(180:\ra) --++(300:\ra) -- cycle;
		\path (240:\ra) --++(300:\ra) coordinate (d1);
		\foreach \w in {1,6}{
			\fill[black!80] (d1) --++(60+60*\w:\ra) --++(180+60*\w:\ra) -- cycle;
		}
		\foreach \w in {1,6}{
			\draw[black!25] (d1) --++(60+60*\w:\ra) --++(180+60*\w:\ra) -- (d1);
		}
		\foreach \w in {1,3,4,5}{
			\draw[black!25] (0,0) --++(60*\w:\ra) --++(120+60*\w:\ra) -- (0,0);
		}
		\draw[black!25] (240:\ra) --++(180:\ra) --++(300:\ra) -- (240:\ra);
		\foreach \w in {1,2}{
			\draw[black!25] (180:\ra) --++(60+60*\w:\ra) --++(180+60*\w:\ra) -- (180:\ra);
		}

        \end{tikzpicture}
        \subcaption{3 holes, 19 tiles}
    \end{subfigure}
    \begin{subfigure}[t]{.3\textwidth}
        \centering
        \begin{tikzpicture}
		\foreach \w in {3,...,6}{
		\LSide{\w}{2}
		}
		\foreach \w in {2,3,4,6}{
			\fill[black!80] (0,0) --++(60*\w:\ra) --++(120+60*\w:\ra) -- cycle;
		}
		\foreach \w in {1,2}{
			\fill[black!80] (180:\ra) --++(60+60*\w:\ra) --++(180+60*\w:\ra) -- cycle;
		}
		\foreach \w in {1,2}{
			\fill[black!80] (180+60*\w:\ra) --++(60+120*\w:\ra) --++(180+120*\w:\ra) -- cycle;
		}
		\path (240:\ra) --++(300:\ra) coordinate (d1);
		\foreach \w in {1,...,5}{
			\fill[black!80] (d1) --++(60+60*\w:\ra) --++(180+60*\w:\ra) -- cycle;
		}
		\foreach \w in {1,...,5}{
			\draw[black!25] (d1) --++(60+60*\w:\ra) --++(180+60*\w:\ra) -- (d1);
		}
		\foreach \w in {2,3,4,6}{
			\draw[black!25] (0,0) --++(60*\w:\ra) --++(120+60*\w:\ra) -- (0,0);
		}
		\foreach \w in {1,2}{
			\draw[black!25] (180+60*\w:\ra) --++(60+120*\w:\ra) --++(180+120*\w:\ra) -- (180+60*\w:\ra);
		}
		\foreach \w in {1,2}{
			\draw[black!25] (180:\ra) --++(60+60*\w:\ra) --++(180+60*\w:\ra) -- (180:\ra);
		}

        \end{tikzpicture}
        \subcaption{4 holes, 23 tiles}
    \end{subfigure}
    \begin{subfigure}[t]{.3\textwidth}
        \centering
        \begin{tikzpicture}
		\LSide{1}{3}
		\foreach \w in {1,2,3,5}{
			\fill[black!80] (0,0) --++(60*\w:\ra) --++(120+60*\w:\ra) -- cycle;
		}
		\foreach \w in{1,2,4,5,6}{
		\LTeeth{\w}{2}
		}
		\fill[black!80] (0:\ra) --++(240:\ra) --++(0:\ra) -- cycle;
		\fill[black!80] (300:\ra) --++(240:\ra) --++(0:\ra) -- cycle;
		\path (180:\ra) --++(240:\ra) coordinate (d1);
		\path (180:\ra) --++(120:\ra) coordinate (d2);
		\foreach \w in {1,2,3}{
			\fill[black!80] (d1) --++(60+60*\w:\ra) --++(180+60*\w:\ra) -- cycle;
		}
		\foreach \w in {1,3}{
		\fill[black!80] (d2) --++(180+60*\w:\ra) --++(300+60*\w:\ra) -- cycle;
		}
		\foreach \w in {1,2,3}{
			\fill[black!80] (300+\w*60:\ra) --++(0:\ra) --++(120:\ra) -- cycle;
		}
		\draw[black!25] (0:\ra) --++(240:\ra) --++(0:\ra) -- (0:\ra);
		\foreach \w in {1,2,3}{
			\draw[black!25] (d1) --++(60+60*\w:\ra) --++(180+60*\w:\ra) -- (d1);
		}
		\foreach \w in {1,3}{
		\draw[black!25] (d2) --++(180+60*\w:\ra) --++(300+60*\w:\ra) -- (d2);
		}
		\foreach \w in {1,2,3,5}{
			\draw[black!25] (0,0) --++(60*\w:\ra) --++(120+60*\w:\ra) -- (0,0);
		}
		\draw[black!25] (300:\ra) --++(240:\ra) --++(0:\ra) -- (300:\ra);
		\foreach \w in {1,2,3}{
			\draw[black!25] (300+\w*60:\ra) --++(0:\ra) --++(120:\ra) --++(240:\ra);
		}
        \end{tikzpicture}
        \subcaption{5 holes, 27 tiles}
    \end{subfigure}

\vspace{1cm}

    \begin{subfigure}[t]{.3\textwidth}
        \centering
        \begin{tikzpicture}
		\foreach \w in {1,2}{
			\LSide{\w}{3}
		}
		\foreach \w in {1,2,3,5}{
			\fill[black!80] (0,0) --++(60*\w:\ra) --++(120+60*\w:\ra) -- cycle;
		}
		\LKTeeth{2}
		\fill[black!80] (0:\ra) --++(240:\ra) --++(0:\ra) -- cycle;
		\fill[black!80] (300:\ra) --++(240:\ra) --++(0:\ra) -- cycle;
		\path (180:\ra) --++(120:\ra) coordinate (d1);
		\foreach \w in {1,2,3}{
			\fill[black!80] (d1) --++(60*\w:\ra) --++(120+60*\w:\ra) -- cycle;
		}
		\foreach \w in {1,2,3}{
			\fill[black!80] (300+\w*60:\ra) --++(0:\ra) --++(120:\ra) -- cycle;
		}
		\draw[black!25] (0:\ra) --++(240:\ra) --++(0:\ra) -- (0:\ra);
		\foreach \w in {1,2,3}{
			\draw[black!25] (d1) --++(60*\w:\ra) --++(120+60*\w:\ra) -- (d1);
		}
		\foreach \w in {1,2,3,5}{
			\draw[black!25] (0,0) --++(60*\w:\ra) --++(120+60*\w:\ra) -- (0,0);
		}
		\draw[black!25] (300:\ra) --++(240:\ra) --++(0:\ra) -- (300:\ra);
		\foreach \w in {1,2,3}{
			\draw[black!25] (300+\w*60:\ra) --++(0:\ra) --++(120:\ra) --++(240:\ra);
		}
        \end{tikzpicture}
        \subcaption{6 holes, 31 tiles}
    \end{subfigure}
    \begin{subfigure}[t]{.3\textwidth}
        \centering
        \begin{tikzpicture}
		\foreach \w in {1,2,3}{
			\LSide{\w}{3}
		}
		\foreach \w in {1,...,5}{
			\fill[black!80] (0,0) --++(300+60*\w:\ra) --++(60+60*\w:\ra) -- cycle;
		}
		\path (60:\ra) --++(120:\ra) coordinate (d1);
		\foreach \w in {1,2,3}{
			\fill[black!80] (d1) --++(300+60*\w:\ra) --++(60+60*\w:\ra) -- cycle;
		}
		\foreach \w in {1,2}{
			\fill[black!80] (300+\w*60:\ra) --++(0:\ra) --++(120:\ra) -- cycle;
		}
		\LKTeeth{2}
		\fill[black!80] (0:\ra) --++(240:\ra) --++(0:\ra) -- cycle;
		\draw[black!25] (0:\ra) --++(240:\ra) --++(0:\ra) -- (0:\ra);
		\foreach \w in {1,...,5}{
			\draw[black!25] (0,0) --++(300+60*\w:\ra) --++(60+60*\w:\ra) -- (0,0);
		}
		\foreach \w in {1,2,3}{
			\draw[black!25] (d1) --++(300+60*\w:\ra) --++(60+60*\w:\ra) -- (d1);
		}
		\foreach \w in {1,2}{
			\draw[black!25] (300+\w*60:\ra) --++(0:\ra) --++(120:\ra) --++(240:\ra);
		}

        \end{tikzpicture}
        \subcaption{7 holes, 35 tiles}
    \end{subfigure} 
    \begin{subfigure}[t]{.3\textwidth}
        \centering
        \begin{tikzpicture}
		\foreach \w in {1,...,4}{
			\LSide{\w}{3}
		}
		\foreach \w in {1,...,5}{
			\fill[black!80] (0,0) --++(300+60*\w:\ra) --++(60+60*\w:\ra) -- cycle;
		}
		\LKTeeth{2}
		\fill[black!80] (0:\ra) --++(240:\ra) --++(0:\ra) -- cycle;
		\path (0:\ra) --++(60:\ra) coordinate (d1);
		\foreach \w in {1,...,4}{
			\fill[black!80] (d1) --++(180+60*\w:\ra) --++(300+60*\w:\ra) -- cycle;
		}
		\draw[black!25] (0:\ra) --++(240:\ra) --++(0:\ra) -- (0:\ra);
		\foreach \w in {1,...,4}{
			\draw[black!25] (d1) --++(180+60*\w:\ra) --++(300+60*\w:\ra) -- (d1);
		}
		\foreach \w in {1,...,5}{
			\draw[black!25] (0,0) --++(300+60*\w:\ra) --++(60+60*\w:\ra) -- (0,0);
		}

        \end{tikzpicture}
        \subcaption{8 holes, 39 tiles}
    \end{subfigure}
    \caption{Crystallized polyiamonds for $3\le h\le 8$.}
\label{less8}
\end{figure}
\begin{definition}
For $h\ge h_k$, let $T_h \ast Spir_k$ denote the polyiamond in which $T_h$ exchanges its central configuration for $Spir_k$. If $h=h_k$, then $T_{h_k} \ast Spir_k=Spir_k$.
\end{definition}

\begin{figure}[H]
    \centering
    \begin{subfigure}[t]{.4\textwidth}
        \centering
        \begin{tikzpicture}
        \SpirTwo
        \FractalLayers{3}{6}
        \end{tikzpicture}
        \subcaption{$T_{45}$}
    \end{subfigure}
\hspace{1cm}
    \begin{subfigure}[t]{.4\textwidth}
        \centering
        \begin{tikzpicture}
        \GSpirK{4}
        \FractalLayers{5}{6}
        \end{tikzpicture}
        \subcaption{$T_{45} \ast Spir_{4}$}
    \end{subfigure}
\caption{The center of $T_{45}$ is exchanged with $Spir_4$.}
\label{spiralswap}
\end{figure}

\begin{corollary}
For $k\ge3$ and $h\ge h_k$, there are at least $k-1$ distinct crystallized polyiamonds with $h$ holes, given by the set $$\{T_h\}\cup\{T_h \ast Spir_l:3\le l\le k\}.$$
\end{corollary}
\noindent This set is depicted in full for $h=h_7=63$ in Figure \ref{63spirals}, and in Figure \ref{bigcrsytal} we take $T_{315}$ from Figure \ref{bigfig} and form $T_{315}\ast Spir_8$.

\begin{figure}[H]
    \centering
        \begin{subfigure}[t]{.45\linewidth}   
        \centering
        \begin{tikzpicture}
            \SpirTwo
            \FractalLayers{3}{7}
        \end{tikzpicture}
    \caption{$T_{63}$}
    \end{subfigure}
\hfill
    \begin{subfigure}[t]{.45\linewidth}    
        \centering
        \begin{tikzpicture}
            \SpirThree
            \FractalLayers{4}{7}
        \end{tikzpicture}
    \caption{$T_{63} \ast Spir_{3}$}
    \end{subfigure}

\vspace{.8cm}

    \begin{subfigure}[t]{.45\linewidth}  
        \centering
        \begin{tikzpicture}
            \SpirK{4}
            \FractalLayers{5}{7}
        \end{tikzpicture}
    \caption{$T_{63} \ast Spir_{4}$}
    \end{subfigure}
\hfill
    \begin{subfigure}[t]{.45\linewidth}    
        \centering
        \begin{tikzpicture}
            \SpirK{5}
            \FractalLayers{6}{7}
        \end{tikzpicture}
    \caption{$T_{63} \ast Spir_{5}$}
    \end{subfigure}

\vspace{.8cm}

    \begin{subfigure}[t]{.45\linewidth}  
        \centering
        \begin{tikzpicture}
            \SpirK{6}
            \FractalLayer{7}
        \end{tikzpicture}
    \caption{$T_{63} \ast Spir_{6}$}
    \end{subfigure}
\hfill
    \begin{subfigure}[t]{.45\linewidth}    
        \centering
        \begin{tikzpicture}
            \SpirK{7}
        \end{tikzpicture}
    \caption{$Spir_{7}$}
    \end{subfigure}
\caption{Six distinct crystallized polyiamonds with 63 holes.}
\label{63spirals}
\end{figure}
\section{Concluding Remarks}
In this paper we have completely solved the problem of finding $g_{\triangle}(h)$ for all $h\geq1$, the minimum number of triangles needed for constructing a polyiamond with $h$ holes. This immediately determines the values of $\ft(n)$ as well, the maximum number of holes that can be enclosed by $n$ tiles. Along with this, we examined several structural conditions of polyiamonds, such as having a dual graph which is a tree, having only holes with an area of 1, and having minimal outer perimeter. We refer to polyiamonds satisfying this trio of conditions as being \textit{efficiently structured}, and have shown here that this is in fact equivalent to the property of being \textit{crystallized}. This, however, is not the case for polyominoes, where being efficiently structured was recently shown to be a stronger condition \cite{malen2019topology}. 

In a more general setting, we can ask these extremal geometrical and topological questions for higher dimensional simplicial complexes, and also consider the effect of changing the dimension of the ambient space.  
\begin{figure}[H]
    \centering
        \begin{tikzpicture}
            \centering
            \SpirK{8}
            \FractalLayers{9}{15}
        \end{tikzpicture}
    \caption{$T_{315}\ast Spir_{8}$}
\label{bigcrsytal}
\end{figure}
\section*{Acknowledgments.}
This work was supported by HFSP RGP0051/2017, NSF DMS 17-13012, and NSF DMS-1352386.
\bibliographystyle{plain}
\bibliography{bibliography}


\end{document}